\documentclass[10pt,reqno]{amsart}

\usepackage{amsmath,amssymb}
\usepackage[all]{xy}
\usepackage{mathbbol}
\usepackage{mathrsfs}

\newtheorem{proposition}{Proposition}[section]
\newtheorem{theorem}[proposition]{Theorem}

\theoremstyle{remark}
\newtheorem{remark}[proposition]{Remark}

\newcommand{\cst}{\ifmmode\mathrm{C}^*\else{$\mathrm{C}^*$}\fi}
\newcommand{\tens}{\otimes}
\newcommand{\atens}{\otimes_{\text{\rm\tiny{alg}}}}

\newcommand{\adelta}{\delta_{\text{\rm\tiny alg}}}
\newcommand{\aDelta}{\Delta_{\text{\rm\tiny alg}}}
\newcommand{\deltau}{\delta_{\text{\rm\tiny{u}}}}
\newcommand{\Deltar}{\Delta_{\text{\rm\tiny{r}}}}
\newcommand{\deltaR}{\delta^{\text{\rm\tiny{r}}}}
\newcommand{\Deltau}{\Delta_{\text{\rm\tiny{u}}}}

\newcommand{\deltaU}{\delta^{\text{\rm\tiny{u}}}}
\newcommand{\piu}{\pi_{\text{\rm\tiny{u}}}}
\newcommand{\Bu}{B_{\text{\rm\tiny{u}}}}
\newcommand{\Au}{A_{\text{\rm\tiny{u}}}}
\newcommand{\Ar}{A_{\text{\rm\tiny{r}}}}
\newcommand{\GGu}{\GG_{\text{\rm\tiny{u}}}}
\newcommand{\deltar}{\delta_{\text{\rm\tiny{r}}}}
\newcommand{\pir}{\pi_{\text{\rm\tiny{r}}}}
\newcommand{\Br}{B_{\text{\rm\tiny{r}}}}
\newcommand{\cBr}{\cB_{\text{\rm\tiny{r}}}}
\newcommand{\GGr}{\GG_{\text{\rm\tiny{r}}}}
\newcommand{\comp}{\circ}

\newcommand{\GG}{\mathbb{G}}

\newcommand{\id}{\mathrm{id}}
\newcommand{\I}{\mathbb{1}}

\newcommand{\cA}{\mathscr{A}}
\newcommand{\cB}{\mathscr{B}}
\newcommand{\tcB}{\widetilde{\mathscr{B}}}
\newcommand{\cR}{\mathscr{R}}
\newcommand{\cS}{\mathscr{S}}
\newcommand{\st}{\:\vline\:}

\DeclareMathOperator{\Mor}{Mor}
\DeclareMathOperator{\spn}{span}

\newcommand{\rB}{\!\sideset{_{\text{\rm\tiny{r}}}}{}{\operatorname{\mathit{B}}}}
\newcommand{\rdelta}{\!\sideset{_{\text{\rm\tiny{r}}}}{}{\operatorname{\delta}}}
\newcommand{\rA}{\!\sideset{_{\text{\rm\tiny{r}}}}{}{\operatorname{\mathit{A}}}}
\newcommand{\rDelta}{\!\sideset{_{\text{\rm\tiny{r}}}}{}{\operatorname{\Delta}}}
\newcommand{\rGG}{\!\sideset{_{\text{\rm\tiny{r}}}}{}{\operatorname{\mathbb{G}}}}

\numberwithin{equation}{section}

\begin{document}

\title{On actions of compact quantum groups}

\date{}

\author{Piotr M.~So{\l}tan}
\address{Institute of Mathematics, Polish Academy of Sciences\newline
and\newline
Department of Mathematical Methods in Physics, Faculty of Physics, University of Warsaw}
\email{piotr.soltan@fuw.edu.pl}

\thanks{Research partially supported by Polish government grant no.~N201 1770 33, European Union grant PIRSES-GA-2008-230836 and Polish government matching grant no.~1261/7.PR UE/2009/7.}

\begin{abstract}
We compare algebraic objects related to a compact quantum group action on a unital $\mathrm{C}^*$-algebra in the sense of Podle\'s and Baum \emph{et al}.~and show that they differ by the kernel of the morphism describing the action. Then we address ways to remove the kernel without changing the Podle\'s algebraic core. A minimal such procedure is described. We end the paper with a natural example of an action of a reduced compact quantum group with non-trivial kernel.
\end{abstract}

\maketitle

\section{Introduction}

In this paper we discuss some aspects of the theory of actions of compact quantum groups on unital \cst-algebras (compact quantum spaces) which are usually neglected or excluded from consideration by additional assumptions. The motivation for this came from some recent work of P.~Baum \emph{et al}.~on the Peter-Weyl functor (\cite{hajac}, cf.~also Section \ref{algCore}).

\subsection{Standing assumptions and notation}

All considered \cst-algebras will be unital and for such \cst-algebras $A$ and $B$ we will denote by $\Mor(A,B)$ the set of unital $*$-homomorphisms from $A$ to $B$. The symbol ``$\tens$'' will denote the minimal tensor product of \cst-algebras.

Throughout the paper $\GG=(A,\Delta)$ will be a compact quantum group as defined in \cite[Definition 2.1]{cqg}. By $\GGu=(\Au,\Deltau)$ and $\GGr=(\Ar,\Deltar)$ we shall denote the \emph{universal} and \emph{reduced} versions of $\GG$ \cite[Sections 3 \& 2]{bmt}. The canonical morphisms from $\Au$ to $A$ and from $A$ to $\Ar$ will be denoted by
\[
\Lambda\in\Mor(\Au,A),\qquad\lambda\in\Mor(A,\Ar).
\]
The canonical dense Hopf $*$-algebra inside $A$ will be denoted by $(\cA,\aDelta)$. In particular we have $\aDelta=\bigl.\Delta\bigr|_{\cA}$ (\cite[Section 6]{cqg}).

Also throughout the paper $B$ will denote a unital \cst-algebra.

\subsection{Actions of compact quantum groups on \cst-algebras}\label{actS}

We shall denote by $\delta$ an \emph{action} of $\GG$ on $B$. This means that $\delta\in\Mor(B,B\tens{A})$ is such that
\begin{enumerate}
\item $(\delta\tens\id)\comp\delta=(\id\tens\Delta)\comp\delta$,
\item $\spn\bigl\{\delta(b)(\I\tens{a})\st{a}\in{A},\;b\in{B}\bigr\}$ is dense in $B\tens{A}$.
\end{enumerate}
Actions of compact quantum groups were defined an studied first by Piotr Podle\'s in \cite{podPhd}. There is now ample literature on this topic, eg.~\cite{boca,wang,li,marciniak,so3}. Some authors assume that $\delta$ is an injective map. Others impose the formally stronger condition that
\begin{equation}\label{epsi}
(\id\tens\epsilon)\comp\delta=\id,
\end{equation}
where $\epsilon$ is the \emph{counit} of $\GG$ (which is then assumed to be continuous). Clearly \eqref{epsi} implies injectivity of $\delta$.

As mentioned in the introduction, we will focus on situations, when such additional conditions do not hold. In fact there are natural examples of actions which are not injective (cf. Section \ref{nonInj}).

\begin{remark}\label{injRem}
\noindent\begin{enumerate}
\item\label{injRem1} 
It is important to note that in case $\GG$ does possess a continuous counit then it follows from our definition of an action that in fact \eqref{epsi} is automatically satisfied. In our case this can be easily proved using the Podle\'s algebraic core (see Section \ref{algCore}), but can be established in a much more general situation (\cite[Lemma 2.2]{skalski}).
\item There are many other situations when an action $\delta\in\Mor(B,B\tens{A})$ is automatically injective. For example this happens if $B$ is finite dimensional or simple (cf.~\cite[Example 3.6]{li}).
\end{enumerate}
\end{remark}

\section{Two approaches to the algebraic core}\label{algCore}

The notion of algebraic core of an action of a compact quantum group goes back to the PhD thesis of Podle\'s (\cite{podPhd}). Motivated by the fact that any action of a compact group on a Banach space decomposes into isotypical components Podle\'s showed that the vector space $\cB$ spanned by elements of $B$ which transform according to irreducible representations of $\GG$ is a dense unital $*$-subalgebra of $B$ which is a right comodule algebra for the Hopf algebra $(\cA,\aDelta)$.

We shall now explain briefly the construction of $\cB$ and introduce notation needed in what follows. Let $\cR$ be a set indexing the equivalence classes of irreducible unitary representations of $\GG$ (cf.~\cite{cqg}). For each $\alpha\in\cR$ we let
\[
u^\alpha=\begin{bmatrix}u^\alpha_{1,1}&\cdots&u^\alpha_{1,N_\alpha}\\
\vdots&\ddots&\vdots\\u^\alpha_{N_\alpha,1}&\cdots&u^\alpha_{N_\alpha,N_\alpha}
\end{bmatrix}\in{M_{N_\alpha}\tens{A}}
\]
be a representative of the class corresponding to $\alpha$. Note that
\begin{equation}\label{bSLW}
\bigl\{u^\alpha_{i,j}\st\alpha\in\cR,\;i,j=1,\ldots,N_\alpha\bigr\}
\end{equation}
is a linear basis of $\cA$ (\cite[Proposition 6.1]{cqg}). There are (continuous) functionals $\phi^\alpha_{i,j}$ on $A$ such that
\[
\phi^\alpha_{i,j}(u^\beta_{k,l})=\delta_{\alpha,\beta}\delta_{i,k}\delta_{j,l}.
\]
We define $E^\alpha:B\to{B}$ by $E^\alpha=\sum\limits_{r=1}^{N_\alpha}(\id\tens\phi^\alpha_{r,r})\comp\delta$ and let
\[
W_\alpha=E^\alpha(B),\quad\cB=\bigoplus_{\alpha\in\cR}W_\alpha\subset{B}.
\]

The following theorem is due to Podle\'s (\cite{podles,podmu}, see also \cite{boca}):

\begin{theorem}
\noindent\begin{enumerate}
\item $\cB$ is a dense unital $*$-subalgebra of $B$,
\item $\delta(\cB)\subset\cB\atens\cA$,
\item $\adelta=\bigl.\delta\bigr|_\cB$ is a coaction of the Hopf $*$-algebra $\cA$ on the $*$-algebra $\cB$, in particular
\[
(\id\tens\epsilon)\delta(b)=b
\]
($\epsilon$ is the counit of $\cA$) for all $b\in\cB$, so that $\ker{\delta}\cap\cB=\{0\}$.
\end{enumerate}
\end{theorem}

We call the algebra $\cB$ the \emph{Podle\'s subalgebra} of $B$ and the comodule algebra
\[
\adelta:\cB\longrightarrow\cB\tens\cA
\]
the \emph{algebraic core} of the action $\delta$.

In the chapter from an upcoming book \cite{hajac} P.~Baum, P.M.~Hajac, R.~Matthes and W.~Szymanski introduce a different $*$-subalgebra of $B$:
\[
\tcB=\bigl\{b\in{B}\st\delta(b)\in{B}\atens\cA\bigr\}
\]
(\cite[Eq.~(3.1.4)]{hajac}). This space is essential in the study of the Peter-Weyl functor (see \cite{hajac}).

It is fairly obvious that $\cB\subset\tcB$. Our aim is to describe the link between $\cB$ and $\tcB$.

\begin{proposition}
\noindent
\begin{enumerate}
\item\label{jeden} For any $b\in\tcB$ we have $\delta(b)\in\cB\atens\cA$,
\item\label{dwa} $\tcB=\cB\oplus\ker{\delta}$.
\end{enumerate}
\end{proposition}

\begin{proof}
Ad \eqref{jeden}. Using the basis \eqref{bSLW} we can write the element $\delta(b)\in{B}\atens\cA$ in the form 
\[
\delta(b)=\sum_{\alpha\in\cS}\sum_{i,j=1}^{N_\alpha}b^\alpha_{i,j}\tens{u^\alpha_{i,j}},
\]
where $\cS$ is a finite subset of $\cR$. Since $(\delta\tens\id)\comp\delta=(\id\tens\Delta)\comp\delta$, we have
\begin{equation}\label{sumsum}
\sum_{\alpha\in\cS}\sum_{i,j=1}^{N_\alpha}\delta(b^\alpha_{i,j})\tens{u^\alpha_{i,j}}
=\sum_{\alpha\in\cS}\sum_{i,j,s=1}^{N_\alpha}b^\alpha_{i,j}\tens{u^\alpha_{i,s}}\tens{u^\alpha_{s,j}}.
\end{equation}
Applying $(\id\tens\id\tens\phi^\beta_{k,l})$ to both sides of \eqref{sumsum} gives
\begin{equation}\label{deb}
\delta(b^\beta_{k,l})=\sum_{i=1}^{N_\beta}b^\beta_{i,l}\tens{u^\beta_{i,k}}
\end{equation}
and it follows that $E(b^\beta_{k,l})=b^\beta_{k,l}$, so that $b^\beta_{k,l}\in{W_\beta}$. This means that $\delta(b)\in\cB\atens\cA$.

Ad \eqref{dwa}. Take $b\in\tcB$ and write
\[
\delta(b)=\sum_{\alpha\in\cS}\sum_{i,j=1}^{N_\alpha}b^\alpha_{i,j}\tens{u^\alpha_{i,j}},
\]
as in the proof of Statement \eqref{jeden}. Let
\[
b'=\sum_{\alpha\in\cS}\sum_{i=1}^{N_\alpha}b^\alpha_{i,i}.
\]
By \eqref{deb}
\[
\delta(b')=\sum_{\alpha\in\cS}\sum_{i=1}^{N_\alpha}\sum_{k=1}^{N_\alpha}b^\alpha_{k,i}\tens{u^\alpha_{k,i}}.
\]
It follows that $\delta(b)=\delta(b')$, so that $b-b'\in\ker{\delta}$. Moreover, by Statement \eqref{jeden}, we have
$b'\in\cB$.
\end{proof}

It is desirable in some applications that $\tcB=\cB$. Section \ref{minRed} is devoted to possible ways of obtaining this equality.

Let us end this section with the following remark:

\begin{remark}
We know that given an action $\delta\in\Mor(B,B\tens{A})$ as defined in Subsection \ref{actS} the algebra $\tcB=\delta^{-1}(B\atens\cA)$ is dense in $B$. Assume that $\GG$ has a continuous counit and consider a map $\theta\in\Mor(B,B\tens{A})$ satisfying only $(\theta\tens\id)\comp\theta=(\id\tens\Delta)\comp\theta$. If $\theta^{-1}(B\atens\cA)$ is dense in $B$ then the condition
\[
(\id\tens\epsilon)\comp\theta=\id
\]
implies that the linear span of $\bigl\{\theta(b)(\I\tens{a})\st{a}\in{A},\;b\in{B}\bigr\}$ is dense in $B\tens{A}$. Indeed, take $c\in\theta^{-1}(B\atens\cA)$ and $a\in\cA$. Let $S$ be the antipode of $\cA$. Then the element
\[
\bigl((\id\tens{S})\delta(c)\bigr)(\I\tens{a})\in{B}\atens\cA
\]
and applying to it the map $\Phi:B\atens{A}\ni(x\tens{y})\mapsto\delta(x)(\I\tens{y})\in{B}\tens{A}$ we get $\Phi(X)=c\tens{a}$. It follows that the range of $\Phi$ is dense in $B\tens{A}$. Note that if $B$ is generated by elements of some set $\mathscr{S}\subset{B}$ then it is enough to check that $\theta(s)\in{B}\atens\cA$ for all $s\in\mathscr{S}$ to ensure that $\theta^{-1}(B\atens\cA)$ is dense in $B$.
\end{remark}

\section{Universal and reduced action. Minimal reduction}\label{minRed}

Given the action $\delta\in\Mor(B,B\tens{A})$ one can perform certain operations on $B$ which are know as passage to the \emph{universal} of \emph{full} action and \emph{reduction} of the action respectively. These have been known for quite some time (cf.~eg.~\cite[Section 1]{boca}). A very good descriptions of both operations can be found in \cite[Section 3]{li}.

We have

\begin{theorem}[{\cite[Proposition 3.3]{li}}]
\noindent\begin{enumerate}
\item The $*$-algebra $\cB$ admits the universal enveloping \cst-algebra $\Bu$;
\item the natural extension of $\adelta$ to a map $\deltau\in\Mor(\Bu,\Bu\tens{A})$ is an action of $\GG$ on $\Bu$;
\item the canonical morphism $\piu\in\Mor(\Bu,B)$ is $\GG$-equivariant, i.e.~the diagram
\[
\xymatrix{
\Bu\ar[rr]^-{\deltau}\ar[d]_{\piu}&&\Bu\tens{A}\ar[d]^{\piu\tens\id}\\
B\ar[rr]_-\delta&&B\tens{A}\\
}
\]
is commutative;
\item the Podle\'s algebra of $\deltau$ is $\cB\subset\Bu$ and the algebraic part of $\deltau$ can be canonically identified with $\adelta$.
\end{enumerate}
\end{theorem}

Let s remark that the passage to the action $\deltau$ we do not get rid of the kernel of the original action $\delta$. Indeed if we take $\GG=\GGu$ in the example in Section \ref{nonInj} then the constructed action is universal, but has a non-zero kernel. On the other hand, the procedure of reduction described in the next theorem leads to an injective action:

\begin{theorem}[{\cite[Proposition 3.4]{li}}]
There exists a unital \cst-algebra $\Br$ and a surjective $*$-homomorphism $\pir\in\Mor(B,\Br)$ such that
\noindent\begin{enumerate}
\item there exists a unique $\deltar\in\Mor(\Br,\Br\tens{A})$ such that
\begin{equation}\label{pir}
\xymatrix{
B\ar[rr]^-{\delta}\ar[d]_{\pir}&&B\tens{A}\ar[d]^{\pir\tens\id}\\
\Br\ar[rr]_-\deltar&&B\tens{A}\\
}
\end{equation}
and $\deltar$ is an action of $\GG$ on $B_r$,
\item $\ker{\deltar}=\{0\}$,
\item $\pir$ is injective on $\cB$,
\item the Podle\'s algebra $\cBr$ of $\deltar$ is equal to $\pir(\cB)$,
\item the algebraic part of $\deltar$ is can be canonically identified with $\adelta$.
\end{enumerate}
\end{theorem}

The procedure of reduction is therefore one way to ensure the equality $\tcB=\cB$. However it can easily happen that for injective $\delta$ (so in a case when we already have $\tcB=\cB$) the procedure of reduction changes $B$ which is not necessary. Indeed, suppose $\GG$ is not reduced ad has a continuous counit. If we take $B=A$ and $\delta=\Delta$ then we have $\Br=\Ar$ (\cite[Example 3.6(3)]{li}) which is a proper quotient of $B$, while the action $\delta$ is injective by Remark \ref{injRem}\eqref{injRem1}.

Below we address the procedure of \emph{minimal reduction}. Let $\rB=B/\ker{\delta}$ and let $p:B\to\rB$ be the quotient map. The $*$-homomorphism $\delta:B\to{B\tens{A}}$ is a composition
\[
\xymatrix{
B\ar[rr]^-\delta\ar[dr]_p&&B\tens{A}\\
&\rB\ar[ur]_{i}
}
\]
with $i$ injective. Define $\rdelta=(p\tens\id)\comp{i}$.

\begin{theorem}
\noindent
\begin{enumerate}
\item\label{pierwsze} $\rdelta$ is an action of $\GG$ on $\rB$ and the morphism $p\in\Mor(B,\rB)$ is equivariant,
\item\label{drugie} $p$ is injective on $\cB$ and the algebraic part of $\rdelta$ can be canonically identified with $\adelta$,
\item\label{trzecie} if $\ker{\Delta}=\{0\}$ then $\ker{\rdelta}=\{0\}$.
\end{enumerate}
\end{theorem}

\begin{proof}
Ad \eqref{pierwsze}. Take $b\in{B}$
\[
\begin{split}
(\rdelta\tens\id)\rdelta\bigl(p(b)\bigr)&=(p\tens\id\tens\id)(i\tens\id)(p\tens\id)i\bigl(p(b)\bigr)\\
&=(p\tens\id\tens\id)(\delta\tens\id)\delta(b)\\
&=(p\tens\id\tens\id)(\id\tens\Delta)\delta(b)\\
&=(p\tens\id\tens\id)(\id\tens\Delta)i\bigl(p(b)\bigr)\\
&=(\id\tens\Delta)(p\tens\id)i\bigl(p(b)\bigr)\\
&=(\id\tens\Delta)\rdelta\bigl(p(b)\bigr).
\end{split}
\]
Therefore we have $(\rdelta\tens\id)\comp\rdelta=(\id\tens\Delta)\comp\rdelta$.

Now since
\[
\begin{split}
\spn\bigl\{\rdelta(x)(\I\tens{a})\st{a}\in{A},\;x\in\rB\bigr\}
&=\spn\bigl\{\bigl((p\tens\id)i(x)\bigr)(\I\tens{a})\st{a}\in{A},\;x\in\rB\bigr\}\\
&=\spn\bigl\{(p\tens\id)\bigl(i(x)(\I\tens{a})\bigr)\st{a}\in{A},\;x\in\rB\bigr\}\\
&=\spn\bigl\{(p\tens\id)\bigl[i\bigl(p(b)\bigr)(\I\tens{a})\bigr]\st{a}\in{A},\;b\in{B}\bigr\}\\
&=\spn\bigl\{(p\tens\id)\bigl(\delta(b)(\I\tens{a})\bigr)\st{a}\in{A},\;b\in{B}\bigr\}\\
&=(p\tens\id)\bigl(\spn\bigl\{\delta(b)(\I\tens{a})\st{a}\in{A},\;b\in{B}\bigr\}\bigr)
\end{split}
\]
which is dense in $\rB\tens{A}$ because $(p\tens\id)$ is surjective. This establishes that $\rdelta$ is an action of $\GG$ on $\rB$.

We have the commutative diagram:
\begin{equation}\label{ipd}
\xymatrix{
B\ar[d]_p\ar[rr]^-\delta&&B\tens{A}\ar[d]^{p\tens\id}\\
\rB\ar[urr]^-i\ar[rr]_-{\rdelta}&&\rB\tens{A}
}
\end{equation}
(the lower triangle is the definition of $\rdelta$) which shows that $p$ is equivariant.

Ad \eqref{drugie}. Since $\cB\cap\ker{\delta}=\{0\}$ we have $\ker{\bigl.p\bigr|_\cB}=\{0\}$. The second assertion follows from the diagram \eqref{ipd}.

Ad \eqref{trzecie}. Let us first note that $\ker{p}=\ker{\delta}$ implies that $\ker(p\tens\id)=\ker(\delta\tens\id)$. To see this note that $(\delta\tens\id)=(i\tens\id)(p\tens\id)$ and $(i\tens\id)$ is injective (\cite[Proposition 4.22]{tak1}).
\[
\begin{split}
\ker{\rdelta}&=\bigl\{x\in\rB\st(p\tens\id)i(x)=0\bigr\}\\
&=\bigl\{p(b)\st{b}\in{B},\;(p\tens\id)i\bigl(p(b)\bigr)=0\bigr\}\\
&=\bigl\{p(b)\st{b}\in{B},\;(p\tens\id)\delta(b)=0\bigr\}\\
&=\bigl\{p(b)\st{b}\in{B},\;\delta(b)\in\ker{(p\tens\id)}\bigr\}\\
&=\bigl\{p(b)\st{b}\in{B},\;\delta(b)\in\ker{(\delta\tens\id)}\bigr\}\\
&=\bigl\{p(b)\st{b}\in{B},\;(\delta\tens\id)\delta(b)=0\bigr\}\\
&=\bigl\{p(b)\st{b}\in{B},\;(\id\tens\Delta)\delta(b)=0\bigr\}\\
&=\bigl\{p(b)\st{b}\in{B},\;\delta(b)=0\bigr\}=\{0\}.
\end{split}
\]
\end{proof}

We now see that if $\ker{\Delta}=\{0\}$ then the minimal reduction is the most economical way to obtain equality $\tcB=\cB$.

There are no known examples of compact quantum groups with non-injective coproduct and we could venture a conjecture that such examples do not exist. It is known that the coproduct of reduced and universal quantum groups is always injective.

\begin{remark}
\noindent\begin{enumerate}
\item
There is a canonical surjective morphism $\gamma:\rB\to\Br$. Indeed, if $b\in\ker{\delta}$ then $\deltar\bigl(\pir(b)\bigr)=0$ (by \eqref{pir}). Now since $\ker{\deltar}=\{0\}$, we have that $b\in\ker{\pir}$. In other words $\ker{\delta}\subset\ker{\pir}$. Moreover $\gamma$ is equivariant because it is an isomorphism of the algebraic parts of $\rdelta$ and $\deltar$. Also $\gamma\comp{p}=\pir$.
\item If $\ker{\Delta}\neq\{0\}$ then one can consider the algebra $\rA=A/\ker{\Delta}$. It is easy to see that we obtain a comultiplication $\rDelta\in\Mor(\rA,\rA\tens\rA)$ and that $\rGG=(\rA,\rDelta)$ is a compact quantum group whose quotient is $\GGr$. However, it is not clear if $\ker{\rDelta}$ is different from $\{0\}$.
\end{enumerate}
\end{remark}

\section{Lifts, and restrictions of actions. Construction of non-injective actions}\label{nonInj}

The reduced version $\GGr$ of $\GG$ is a \emph{quantum subgroup} of $\GG$ in the sense that the reduction map $\lambda\in\Mor(A,\Ar)$ intertwined the comultiplications: $\Deltar\comp\lambda=(\lambda\tens\lambda)\comp\Delta$. Clearly an action $\delta\in\Mor(B,B\tens{A})$ can always be \emph{restricted} to a quantum subgroup and thus, in particular, we obtain the restriction of $\delta$ to $\deltaR=(\id\tens\lambda)\comp\delta\in\Mor(B,B\tens\Ar)$ which is an action of $\GGr$ on $B$. The algebraic core of $\deltaR$ can be identified with that of $\delta$.

One of the consequences of this fact is that without changing the algebraic core one can always ensure that a given action $\delta\in\Mor(B,B\tens{A})$ is transformed into an action of $\GGr$ whose coproduct is injective, while the algebraic core remains the same.

Another question which one might ask is whether a given action of $\GG$ admits a \emph{lift} to an action of $\GGu$. In other words if $\delta\in\Mor(B,B\tens{A})$ is an action, is there an action $\deltaU\in\Mor(B,B\tens\Au)$ of $\GGu$ on $B$ such that the diagram
\[
\xymatrix{
&&B\tens\Au\ar[d]^{\id\tens\Lambda}\\
B\ar[urr]^-{\deltaU}\ar[rr]_-{\delta}&&B\tens{A}
}
\]
is commutative. It was shown in \cite[Theorem 4.7]{fischer} that if $\GG$ is reduced then any continuous action of $\GG$ with trivial kernel admits a lift.

\begin{proposition}\label{ker}
Assume that $\GG=(A,\Delta)$ is not reduced. Put $B=A$ and let $\delta=(\id\tens\lambda)\comp\Delta$. Then $\delta\in\Mor(B,B\tens\Ar)$ is an action of $\GGr$ on $B$ and $\ker{\delta}=\ker{\lambda}$. In particular $\delta$ is not injective.
\end{proposition}

\begin{proof}
Let $\pi$ be the composition $\xymatrix@1{\Au\ar[r]^\Lambda&A\ar[r]^\lambda&\Ar}$. It is easy to see that the universal lift of the (injective) action of $\GGr$ on itself given by $\Deltar$ can be performed just as well ``on the left leg''. We obtain a map which we will call $\deltaU$ (just as in the discussion preceding our proposition) for which the diagram
\[
\xymatrix{
&&\Au\tens\Ar\ar[d]^{\pi\tens\id}\\
\Ar\ar[urr]^-{\deltaU}\ar[rr]_-{\Deltar}&&\Ar\tens{\Ar}
}
\]
is commutative.

We then have the following commutative diagram
\[
\xymatrix{
A\ar[dd]_\lambda\ar[rr]^-\Delta\ar[rrd]_-\delta&&A\tens{A}\ar[d]^{\id\tens\lambda}\\
&&A\tens\Ar\ar[d]^{\lambda\tens\id}\\
\Ar\ar[urr]^-{\widetilde{\delta}}\ar[rr]_{\Deltar}&&\Ar\tens\Ar
}
\]
where $\widetilde{\delta}=(\Lambda\tens\id)\comp\deltaU$.

Assume now that $x\in\ker{\lambda}$. Then $\delta(x)=\widetilde{\delta}\bigl(\lambda(x)\bigr)=0$, which shows that $\ker{\lambda}\subset\ker{\delta}$. On the other hand, if $\delta(x)=0$ then
\[
(\lambda\tens\lambda)\Delta(x)=(\lambda\tens\id)\delta(x)=0.
\]
Since $(\lambda\tens\lambda)\comp\Delta=\Deltar\comp\lambda$ and $\Deltar$ is injective, we see that $\lambda(x)=0$, i.e.~we have $\ker{\delta}\subset\ker{\lambda}$.
\end{proof}

Let us end with a remark that if $\GG$ is not co-amenable (i.e.~the canonical map $\Au\to\Ar$ is not an isomorphism, cf.~\cite{bmt}) then Proposition \ref{ker} gives an example of a non-injective action of $\GGr$ which admits a lift to an action of $\GGu$: we take $\GG=\GG_u$ and let $B=\Au$ and $\delta=(\id\tens\lambda)\comp\Deltau$. Then by the proposition $\ker\delta=\ker\lambda\neq\{0\}$ and the lift is provided by $\deltaU=\Deltau:B\to{B}\tens\Au$. This shows that injectivity of the action is not necessary for existence of a lift to the universal level.

\end{document}